\documentclass[12pt,reqno]{amsart}
\usepackage{amsmath, amsfonts, amssymb, amsthm, hyperref}
\usepackage{bm}
\allowdisplaybreaks[4]
\def\l{\left}
\def\r{\right}
\def\bg{\bigg}
\def\({\bg(}
\def\){\bg)}
\def\t{\text}
\def\f{\frac}

\def\ls{\leqslant}

\def\bi{\binom}
\def\al{\alpha}

\def\Proof{\noindent{\it Proof}}

\def\Z{\mathbb Z}
\def\C{\mathbb C}
\def\N{\mathbb N}

\def\1{{\bf 1}}

\def\<{\langle}
\def\>{\rangle}
\theoremstyle{plain}
\newtheorem{theorem}{Theorem}[section]
\newtheorem{lemma}{Lemma}

\theoremstyle{definition}

\theoremstyle{remark}
\newtheorem{remark}{Remark}

\makeatletter
\@namedef{subjclassname@2020}{%
	\textup{2020} Mathematics Subject Classification}
\makeatother
\vspace{4mm}

\begin{document}
	\hbox{Preprint}
	\medskip
	
	\title[On series identities involving $\bi{4k}k$ and harmonic numbers]
	{On series identities involving $\bi{4k}k$ and harmonic numbers}
	\author{Bo Jiang}
	\address {(Bo Jiang) School of Mathematics, Nanjing
		University, Nanjing 210093, People's Republic of China}
	\email{bjiang@smail.nju.edu.cn}
	\author{Zhi-Wei Sun}
	\address {(Zhi-Wei Sun, corresponding author) School of Mathematics, Nanjing
		University, Nanjing 210093, People's Republic of China}
	\email{zwsun@nju.edu.cn}
	
	\keywords{Binomial coefficient, series identity, generating function, harmonic number.
		\newline \indent Supported by the National Natural Science Foundation of China (grant no. 12371004).}
	\begin{abstract}
		The harmonic numbers are those $H_n=\sum_{0<k\ls n}\f1k\ (n=0,1,2,\ldots)$.
		In this paper we confirm over ten conjectural series identities with summands involving the binomial coefficient $\bi{4k}k$ and harmonic numbers. For example, we prove the identities
		\begin{equation*}
			\sum_{k=1}^\infty \frac{\binom{4k}{k}}{16^k}\l((22k^2-92k+11)H_{4k}-\frac{449k-275}{2}-\frac{85}{12k}\r)
			=-151-\frac{80}{3}\log{2}
		\end{equation*}
		and
		\begin{equation*}	 \sum_{k=0}^\infty\frac{\binom{4k}{k}((11k^2+8k+1)(10H_{4k}-17H_{2k})+2k+18)}{(3k+1)(3k+2)16^k}=8\log2,
		\end{equation*}
		which were previously conjectured by Z.-W. Sun.
	\end{abstract}
	\maketitle
	
	\section{Introduction}
	\setcounter{lemma}{0}
	\setcounter{theorem}{0}
	\setcounter{equation}{0}
	\setcounter{conjecture}{0}
	\setcounter{remark}{0}
	\setcounter{corollary}{0}
	
	By Stirling's formula, $n!\sim \sqrt{2\pi n}(n/e)^n$ as $n\to+\infty$.
	Thus, for any integer $m>1$ we have
	$$\bi{mk}k\sim\f{\sqrt m}{\sqrt{2\pi(m-1)k}}\l(\f{m^m}{(m-1)^{m-1}}\r)^k$$
	as $k\to+\infty$. In particular,
	$$\lim_{k\to+\infty}\root k\of{\bi{4k}k}=\f{256}{27}.$$
	
	Z.-W. Sun \cite{1} studied series of the type
	$$\sum_{k=1}^\infty\f{ak^2+bk+c}{k(3k-1)(3k-2)m^k\bi{4k}k},$$
	where $a,b,c$ and $m\not=0$ are rational numbers. For example, Sun \cite{1} deduced the identities
	$$\sum_{k=1}^\infty\frac{(5k^2-4k+1)8^{k}}{k(3k-1)(3k-2)\binom{4k}k}=\frac{3}2\pi$$
	and $$\sum_{k=1}^\infty\frac{415k^2-343k+62}{k(3k-1)(3k-2)(-8)^k\binom{4k}k}=-3\log2$$
	via integration.
	
	Sun's
	conjectural identity
	\begin{equation}\label{-5}
		\sum_{k=0}^\infty\frac{(22k^2-92k+11)\binom{4k}{k}}{16^k}=-5
	\end{equation}
	was confirmed by Max Alekseyev (cf. \cite{2}) via the generating function method.
	Sun \cite{1} proved further that
	\begin{align}
		\label{17}\sum_{k=0}^\infty\f{(22k^2+17k-2)\bi{4k}k}{(k+1)16^k}&=17,
		\\\label{1}\sum_{k=0}^\infty\f{(11k^2+8k+1)\bi{4k}k}{(3k+1)(3k+2)16^k}&=1,
		\\\label{-1/3}\sum_{k=0}^\infty\f{(22k^2-18k+3)\bi{4k}k}{(2k-1)(4k-1)(4k-3)16^k}&=-\f13.
	\end{align}
	
	Recall that the harmonic numbers are those rational numbers
	$$H_n=\sum_{0<k\ls n}\f1k\ \ (n\in\N=\{0,1,2,\ldots\}).$$
	Sun \cite{SSM,1} posed many conjectures on series whose summands involve both $\bi{4k}k$
	and harmonic numbers.
	In this paper, we aim to prove some of such conjectures.
	
	Our first theorem confirms \cite[Conjecture 4.1(i)]{1}.
	\begin{theorem}\label{1}
		Let $P(k)=22k^2-92k+11$. Then
		\begin{equation}
			\sum_{k=1}^\infty \frac{\binom{4k}{k}}{16^k}\l(P(k)H_{k}-54k+108-\frac{10}{3k}\r)=-\frac{20}{3}\log{2}.
			\label{H_k}
		\end{equation}
		\begin{equation}
			\sum_{k=1}^\infty \frac{\binom{4k}{k}}{16^k}\l(P(k)H_{2k}+287k-115-\frac{25}{6k}\r)=214-\frac{40}{3}\log{2},
			\label{H_2k}
		\end{equation}
		\begin{equation}
			\sum_{k=1}^\infty \frac{\binom{4k}{k}}{16^k}\l(P(k)H_{3k}-296k+178-\frac{25}{3k}\r)=-196-\frac{80}{3}\log{2},
			\label{H_3k}
		\end{equation}
		\begin{equation}
			\sum_{k=1}^\infty \frac{\binom{4k}{k}}{16^k}\l(P(k)H_{4k}-\frac{449k-275}{2}-\frac{85}{12k}\r)=-151-\frac{80}{3}\log{2}.
			\label{H_4k}
		\end{equation}
	\end{theorem}
	\begin{remark} The four identities in Theorem \ref{1} were posed to MathOverflow (cf. \cite{MO}) in Feb. 2025,
		but nobody knew how to prove them.
	\end{remark}
	
	Our second theorem confirms the first three identities of \cite[Conjecture 4.2(i)]{1}.
	
	\begin{theorem}\label{2}
		We have the following identities:
		\begin{equation}	 \sum_{k=1}^\infty\frac{\binom{4k}{k}((11k^2+8k+1)H_k+6k+6+4/(3k))}{(3k+1)(3k+2)16^k}=\frac{4}{3}\log2,
			\label{H_k,16,(3k+1)(3k+2)}
		\end{equation}
		\begin{equation}		 \sum_{k=0}^\infty\frac{\binom{4k}{k}((11k^2+8k+1)(H_{2k}-\frac{5}{4}H_k)+4k+1)}{(3k+1)(3k+2)16^k}=\log2,
			\label{4H_2k-5H_k,16,(3k+1)(3k+2)}
		\end{equation}
		and
		\begin{equation}	 \sum_{k=0}^\infty\frac{\binom{4k}{k}((11k^2+8k+1)(10H_{4k}-17H_{2k})+2k+18)}{(3k+1)(3k+2)16^k}=8\log2.
			\label{10H_4k-17H_2k,16,(3k+1)(3k+2)}
		\end{equation}
	\end{theorem}
	\begin{remark} In the spirit of the proof of \cite[Theorem 1.3]{1}, our Theorem \ref{2} actually implies (3.3) and (3.4) in \cite[Conjecture 3.1]{SSM}.
	\end{remark}
	
	Our third theorem confirms \cite[Conjectures 4.3-4.7]{1}.
	
	\begin{theorem}\label{3}
		Let $H(k)=2H_{4k}-3H_{2k}+H_k$	for $k\in\N$. Then we have the following identities:
		\begin{equation}
			 \sum_{k=0}^\infty\frac{\binom{4k}{k}((224k^2-86k+1)H(k)+182k+5)}{(-256)^k}=\frac{9-5\log2}{4\sqrt{2}},
			\label{-256,h}
		\end{equation}
		\begin{equation}
			\sum_{k=0}^\infty \frac{\binom{4k}{k}((112k^2+110k+23)H(k)+28k+16)}{(3k+1)(3k+2)(-256)^k}=8\sqrt{2}\log2,
			\label{-256,3k+2,h}
		\end{equation}
		\begin{equation}
			 \sum_{k=0}^\infty\frac{\binom{4k}{k}((200k^2+76k-17)H(k)-8(725k-49))}{128^k}=\sqrt{2}(144+5\log2),
			\label{128,h}
		\end{equation}
		\begin{equation}
			\sum_{k=0}^\infty \frac{\binom{4k}{k}((40k^2+44k+11)H(k)-8(k+1))}{(3k+1)(3k+2)128^k}=-4\sqrt{2}\log2,
			\label{128,3k+2,h}
		\end{equation}
		\begin{equation}
			 \sum_{k=0}^\infty\frac{\binom{4k}{k}((3575k^2-1026k+67)H(k)+\frac{242}{13}(175k+12))}{(-72)^k}
			=\sqrt{3}\left(\frac{216}{13}-15\log3\right),
			\label{-72,h}
		\end{equation}
		\begin{equation}
			\sum_{k=0}^\infty \frac{\binom{4k}{k}((55k^2+54k+11)H(k)+22k+12)}{(3k+1)(3k+2)(-72)^k}=3\sqrt{3}\log3,
			\label{-72,3k+2,h}
		\end{equation}
		\begin{equation}
			 \sum_{k=0}^\infty\frac{\binom{4k}{k}((21413k^2-1409k+1036)H(k)+\frac{4}{23}(118237k+17320))}{(-25)^k}
			=\sqrt{5}\left(\frac{1440}{23}-100\log 5\right),
			\label{-25,h}
		\end{equation}
		\begin{equation}
			\sum_{k=0}^\infty \frac{\binom{4k}{k}((133k^2+131k+26)H(k)+76k+40)}{(3k+1)(3k+2)(-25)^k}=5\sqrt{5}\log5,
			\label{-25,3k+2,h}
		\end{equation}
		\begin{equation}
			 \sum_{k=0}^\infty\frac{\binom{4k}{k}((49k^2-146k+21)H(k)-3038k+1160)}{24^k}=\sqrt{3}(216-5\log3),
			\label{24,h}
		\end{equation}
		\begin{equation}
			\sum_{k=0}^\infty \frac{\binom{4k}{k}((7k^2+10k+3)H(k)-2k-4)}{(3k+1)(3k+2)24^k}=-\sqrt{3}\log3.
			\label{24,3k+2,h}
		\end{equation}
	\end{theorem}
	
	In the next section we provide some basic lemmas. Theorems \ref{1}-\ref{3}
	will be proved in Sections 3-5 respectively. Our proofs use the basic fact
	\begin{equation}\label{Har}H_n=\sum_{k=1}^n\int_0^1t^{k-1}dt=\int_0^1\sum_{k=1}^n t^{k-1}dt
		=\int_0^1\f{1-t^n}{1-t}dt\ \ \ (n=1,2,\ldots)
	\end{equation}
	and the functional equation of the generating function
	\begin{equation} f(x)=\sum_{k=0}^{\infty}\binom{4k}{k}x^k\quad\l(|x|<\f{27}{256}\r).
	\end{equation}

	\section{Some basic lemmas}
	\setcounter{lemma}{0}
	\setcounter{theorem}{0}
	\setcounter{equation}{0}
	\setcounter{conjecture}{0}
	\setcounter{remark}{0}
	\setcounter{corollary}{0}
	
	\medskip

	\begin{lemma}\label{Lem2.1}
		Let $m\in\Z^+=\{1,2,3,\ldots\}$, and define
		\begin{equation}
			G_m(x)=\sum_{k=0}^{\infty}\frac{1}{(m-1)k+1}\binom{mk}{k}x^k
			\quad \t{for}\ \ |x|<\f{(m-1)^{m-1}}{m^m}. \label{G_m,function}
		\end{equation}
		
		{\rm (i)} $G_m(x)$ satisfies the functional equation
		\begin{equation}
			G_m(x)=1+xG_m(x)^m.\label{G_m}
		\end{equation}

		{\rm (ii)} For $0\ls x<(m-1)^{m-1}/m^m$, we have
		\begin{equation}
			m\log G_m(x)=\sum_{k=1}^{\infty}\frac{1}{k}\binom{mk}{k}x^k.
			\label{logGm(x)}
		\end{equation}
		and
		\begin{equation} \sum_{k=0}^{\infty}\binom{mk}{k}x^k=G_m(x)+(m-1)xG_m'(x)=\frac{G_m(x)}{m-(m-1)G_m(x)}.\label{3.10}
		\end{equation}
	\end{lemma}
	
	\begin{remark} These easy facts are known. For example, parts (i) and (ii) can be found in P. Hilton and J. Pedersen \cite{4}, and Y. Wang, Y. Li and C. Xu \cite[(3.10)-(3.11)]{6}, respectively.
	\end{remark}
	
	\begin{lemma} \label{Lem2.2} Suppose that $|x|<27/256$.
		
		{\rm (i)\ (Max Alekseyev \cite{2})} We have
		\begin{equation}
			27f(x)^4-18f(x)^2-8f(x)-1=256xf(x)^4. \label{equation of f}
		\end{equation}
		
		{\rm (ii)} We have
		\begin{equation}
			\sum_{k=1}^{\infty}\frac{1}{k}\binom{4k}{k}x^k=4\log\frac{4f(x)}{3f(x)+1}.
			\label{f,2}
		\end{equation}
		
		{\rm (iii)} We have
		\begin{equation}
			f'(x)=\frac{64f(x)^5}{(3f(x)+1)^2}
			\label{f,3}
		\end{equation}
		and
		\begin{equation}
			f''(x)=\frac{4096f(x)^9(9f(x)+5)}{(3f(x)+1)^5}.
			\label{f,4}
		\end{equation}
	\end{lemma}
	\Proof.  (i) By (\ref{G_m}) and (\ref{3.10}) with $m=4$, we have
	\begin{equation}
		G_4(x)=1+xG_4(x)^4\ \ \t{and}\ \
		f(x)=\frac{G_4(x)}{4-3G_4(x)}=-\f13+\f4{3(4-3G_4(x))}.\label{f=G_4}
	\end{equation}
	Thus
	\begin{equation}
		\sum_{k=0}^{\infty}\frac{\binom{4k}{k}}{3k+1}x^k=G_4(x)=\frac{4f(x)}{3f(x)+1}
		\label{f,1},
	\end{equation}
	and hence
	\begin{equation}
		\frac{4f(x)}{3f(x)+1}=1+x\l(\frac{4f(x)}{3f(x)+1}\r)^4
	\end{equation}
	which yields (\ref{equation of f}).
	
	{\rm (ii)} (\ref{f,2}) follows from (\ref{logGm(x)}) and
	the second equality in \eqref{f,1}.
	
	{\rm (iii)}
	By taking derivatives with respect to $x$, we obtain from (\ref{equation of f})
	the identities
	\begin{equation}\label{f'}
		f'(x)=\frac{64f(x)^4}{27f(x)^3-9f(x)-2-256xf(x)^3}
	\end{equation}
	and
	\begin{equation}\label{f''} f''(x)=\frac{(-81+768x)f(x)^2f'(x)^2+9f'(x)^2+512f(x)^3f'(x)}{27f(x)^3-9f(x)-2-256xf(x)^3}.
	\end{equation}
	Note that
	\begin{equation}
		27f(x)^3-9f(x)-2-256xf(x)^3=\frac{(3f(x)+1)^2}{f(x)}
	\end{equation}
	by (\ref{equation of f}). Thus we obtain  (\ref{f,3}) and (\ref{f,4})
	from \eqref{f'} and \eqref{f''}.
	
	In view of the above, we have completed the proof of Lemma \ref{Lem2.2}.
	\qed
	
	For any positive integer $j$, we define
	\begin{equation} F_j(x)=\sum_{k=0}^\infty\bi{4k}kH_{jk}x^k\quad\t{for}\ |x|<\f{27}{256}.
	\end{equation}
	
	\begin{lemma}\label{Lem-j} Let $j$ be any positive integer. For $|x|<27/256$, we have
		\begin{equation}
			 F_j(x)=\int_{1}^{f(x)}\frac{4(y-f(x))(1+3y)^2}{jy((3y+1)^3(y-1)-256xy^4\varphi(x,y)^{1-1/j})}dy,
			\label{Fj,d}
		\end{equation}
		\begin{equation}
			 F_j'(x)=\int_{1}^{f(x)}\frac{4((3y^3-2y^2-y)/(4x)-f'(x))(1+3y)^2}{jy((3y+1)^3(y-1)-256xy^4\varphi(x,y)^{1-1/j})}dy,
			\label{Fj,dd}
		\end{equation}
		and
		\begin{equation}
			 F_j''(x)=\int_{1}^{f(x)}\frac{4((27y^5-30y^4-16y^3+14y^2+5y)/(16x^2)-f''(x))(1+3y)^2}{jy((3y+1)^3(y-1)-256xy^4\varphi(x,y)^{1-1/j})}dy,
			\label{Fj,ddd}
		\end{equation}
		where
		\begin{equation} \varphi(x,y)=\frac{(3y+1)^3(y-1)}{256xy^4}.
		\end{equation}
	\end{lemma}
	\Proof. In view of \eqref{Har}, for $|x|<27/256$ we deduce that
	\begin{equation}\label{Fjx} F_j(x)=\sum_{k=1}^\infty\bi{4k}kx^k\int_0^1\f{1-t^{jk}}{1-t}dt
		=\int_0^1\f{f(x)-f(t^jx)}{1-t}dt,
	\end{equation}
	\begin{align*}F_j'(x)&=\sum_{k=1}^\infty\bi{4k}kH_{jk}kx^{k-1}
		=\sum_{k=1}^\infty\bi{4k}kkx^{k-1}\int_0^1\f{1-t^{jk}}{1-t}dt
		 \\&=\int_0^1\f{\sum_{k=1}^\infty\bi{4k}kkx^{k-1}-t^j\sum_{k=1}^\infty\bi{4k}kk(t^jx)^{k-1}}{1-t}dt
	\end{align*}
	and thus
	\begin{equation}\label{Fj'} F_j'(x)=\int_0^1\f{f'(x)-t^jf'(t^jx)}{1-t}dt.
	\end{equation}
	Similarly, for $|x|<27/256$ we have
	\begin{align*}F_j''(x)&=\sum_{k=2}^\infty\bi{4k}kH_{jk}k(k-1)x^{k-2}
		=\sum_{k=2}^\infty\bi{4k}kk(k-1)x^{k-2}\int_0^1\f{1-t^{jk}}{1-t}dt
		 \\&=\int_0^1\f{\sum_{k=2}^\infty\bi{4k}kk(k-1)x^{k-2}-t^{2j}\sum_{k=2}^\infty\bi{4k}kk(k-1)(t^jx)^{k-2}}{1-t}dt
	\end{align*}
	and thus
	\begin{equation}\label{Fj''} F_j''(x)=\int_0^1\f{f''(x)-t^{2j}f''(t^jx)}{1-t}dt.
	\end{equation}
	
	By (\ref{equation of f}), (\ref{f,3}) and \eqref{f,4}, for $0\ls t\ls 1$ we have
	$$
	t^j=\frac{(3f(t^jx)+1)^3(f(t^jx)-1)}{256xf(t^jx)^4}=\varphi(x,f(t^jx)),$$
	and
	$$f'(t^jx)=\frac{64f(t^jx)^5}{(3f(t^jx)+1)^2}
	\ \t{and}\ f''(t^jx)=\frac{4096f(t^jx)^9(9f(t^jx)+5)}{(3f(t^jx)+1)^5}.$$
	If we set $y=f(t^jx)$, then $t^j=\varphi(x,y)$ and $$\frac{dy}{dt}=jf'(t^jx)t^{j-1}x=\frac{64jxy^5}{(3y+1)^2}\varphi(x,y)^{1-1/j}.$$
	Thus, from \eqref{Fjx} we deduce that
	\begin{equation*}
		\begin{aligned}
			 F_j(x)&=\int_{1}^{f(x)}\frac{f(x)-y}{1-\varphi(x,y)^{1/j}}\cdot\frac{(3y+1)^2}{64jxy^5\varphi(x,y)^{1-1/j}}dy\\
			 &=\int_{1}^{f(x)}\frac{4(y-f(x))(1+3y)^2}{jy((3y+1)^3(y-1)-256xy^4\varphi(x,y)^{1-1/j})}dy.
		\end{aligned}
	\end{equation*}
	This proves \eqref{Fj,d}. Similarly, from \eqref{Fj'} we get
	\begin{align*}
		F_j'(x)&=	 \int_{1}^{f(x)}\frac{f'(x)-\varphi(x,y)\frac{64y^5}{(3y+1)^2}}{1-\varphi(x,y)^{1/j}}
		\cdot\frac{(3y+1)^2}{64jxy^5\varphi(x,y)^{1-1/j}}dy\\
		&=\int_{1}^{f(x)}\frac{4((3y^3-2y^2-y)/(4x)-f'(x))(1+3y)^2}{jy((3y+1)^3(y-1)
			-256xy^4\varphi(x,y)^{1-1/j})}dy
	\end{align*}
	and thus (\ref{Fj,dd}) holds. In view of \eqref{Fj''}, we also have
	\begin{equation*}
		\begin{aligned}
			F_j''(x)&= \int_{1}^{f(x)}\frac{f''(x)-\varphi(x,y)^2\frac{4096y^9(9y+5)}{(3y+1)^5}}{1-\varphi(x,y)^{1/j}}
			\cdot\frac{(3y+1)^2}{64jxy^5\varphi(x,y)^{1-1/j}}dy\\
			 &=\int_{1}^{f(x)}\frac{4((27y^5-30y^4-16y^3+14y^2+5y)/(16x^2)-f''(x))(1+3y)^2}{jy((3y+1)^3(y-1)
				-256xy^4\varphi(x,y)^{1-1/j})}dy
		\end{aligned}
	\end{equation*}
	and hence (\ref{Fj,ddd}) is true.
	
	In view of the above, we have completed the proof of Lemma \ref{Lem-j}.
	\qed
	
	\begin{lemma}
		Let
		$$\alpha=f\l(\frac{1}{16}\r),\ \al'=f'\l(\f1{16}\r)
		\ \t{and}\ \al''=f''\l(\f1{16}\r).$$
		Then we have
		\begin{equation}\label{11}
			11\alpha^3-11\alpha^2-7\alpha-1=0
		\end{equation}
		and
		\begin{equation}\label{relation}\frac{11}{128}\al''-\frac{35}{8}\al'+11\al+5=0.
		\end{equation}
	\end{lemma}
	\Proof. Applying \eqref{equation of f} with $x=\frac{1}{16}$,  we get $(11\alpha^3-11\alpha^2-7\alpha-1)(\al+1)=0$.
	Since $\alpha>0$, we see that \eqref{11} holds.
	
	By taking the first derivative and the second derivative of $f(x)$, we see that
	\begin{equation*}
		\frac{11}{128}\al''-\frac{35}{8}\al'+11\al=\sum_{k=0}^\infty \frac{(22k^2-92k+11)\binom{4k}{k}}{16^k}.
	\end{equation*}
	Combining this with \eqref{-5} we obtain \eqref{relation}.
	\qed

	\section{Proof of Theorem $\ref{1}$}
	\setcounter{lemma}{0}
	\setcounter{theorem}{0}
	\setcounter{equation}{0}
	\setcounter{conjecture}{0}
	\setcounter{remark}{0}
	\setcounter{corollary}{0}
	
	\medskip
	
	For convenience, for $j=1,2,3,4$ we set
	$$A_j=F_j\l(\f1{16}\r),\ A_j'=F_j'\l(\f1{16}\r)\ \t{and}\ A_j''=F_j''\l(\f1{16}\r).$$
	
	\medskip
	\noindent{\tt Proof of  $(\ref{H_k})$}.
	Let
	$$\sigma_1:=\sum_{k=1}^\infty \frac{\binom{4k}{k}}{16^k}\l(P(k)H_{k}-54k+108-\frac{10}{3k}\r).$$
	Taking the derivative and the second-order derivative of $F_1(x)$, we deduce that
	\begin{align*}
		\sigma_1
		 =\frac{22}{256}A_1''-\frac{70}{16}A_1'+11A_1-\frac{54}{16}\al'+108(\al-1)-\frac{10}{3}\sum_{k=1}^\infty \frac{\binom{4k}{k}}{k16^k}.
	\end{align*}
	Applying \eqref{Fj,d}--\eqref{Fj,ddd} with $j=1$ and $x=1/{16}$, we get
	\begin{equation} A_1=\int_{1}^{\alpha}\frac{4(y-\alpha)(1+3y)^2}{y(y+1)(11y^3-11y^2-7y-1)}dy,\label{F1,d}
	\end{equation}
	\begin{equation} A_1'=\int_{1}^{\al}\frac{4(4(3y^3-2y^2-y)-\al')(1+3y)^2}{y(y+1)(11y^3-11y^2-7y-1)}dy,
		\label{F1,dd}
	\end{equation}
	and
	\begin{equation} A_1''=\int_{1}^{\al}\frac{4(16(27y^5-30y^4-16y^3+14y^2+5y)-\al'')(1+3y)^2}{y(y+1)(11y^3-11y^2-7y-1)}dy.\label{F1,ddd}
	\end{equation}
	So we deduce further that
	\begin{equation}\label{sigma1}\sigma_1=\int_{1}^{\alpha}\frac{(27y^2-3y-40)(3y+1)^2}{2y(y+1)}dy
		 -\frac{54}{16}\cdot\frac{64\alpha^5}{(3\alpha+1)^2}+108(\alpha-1)-\frac{40}{3}\log\frac{4\alpha}{3\alpha+1}
	\end{equation}
	with the aid of the identity
	\begin{equation}	\begin{aligned}
			&\ \ \ (27y^2-3y-40)(11y^3-11y^2-7y-1)\\
			&=\frac{11}{16}\l(16(27y^5-30y^4-16y^3+14y^2+5y)-\al''\r)
			\\&\ \ \ -35(4(3y^3-2y^2-y)-\al')+11(y-\al).
		\end{aligned}	\label{pqr,1/16}
	\end{equation}
	equivalent to \eqref{relation}.
	For the function
	$$g(y)=\frac{1}{2}(81y^3-54y^2-243y+216)-20\log\frac{2y}{y+1},$$
	it is easy to verify that
	$$g'(y)=\frac{(27y^2-3y-40)(3y+1)^2}{2y(y+1)}.$$
	Thus, from \eqref{sigma1} we obtain
	 \begin{align*}\sigma_1&=g(\al)-g(1)-\frac{216\alpha^5}{(3\alpha+1)^2}+108(\alpha-1)-\frac{40}{3}\log\frac{4\alpha}{3\alpha+1}\\
		&=\frac{27\alpha(\alpha+1)(11\alpha^3-11\alpha^2-7\alpha-1)}{2(3\alpha+1)^2}
		\\&\ \ \ -\frac{20}{3}\log{\l(2+\frac{2(5\alpha^2+2\alpha+1)(11\alpha^3-11\alpha^2-7\alpha-1)}{(\alpha+1)^3(3\alpha+1)^2}\r)} \\
		&=-\frac{20}{3}\log{2}.
	\end{align*}
	This proves (\ref{H_k}).
	\qed

	\noindent{\tt Proof of  $(\ref{H_2k})$}.
	Let
	$$\sigma_2:=\sum_{k=1}^\infty \frac{\binom{4k}{k}}{16^k}\l(P(k)H_{2k}+287k-115-\frac{25}{6k}\r)-214+\frac{40}{3}\log{2}.$$
	Taking the derivative and the second-order derivative of $F_2(x)$, we deduce that
	\begin{align*}
		\sigma_2	 =\frac{22}{256}A_2''-\frac{70}{16}A_2'+11A_2+\frac{287}{16}\al'-115(\al-1)-\frac{25}{6}\sum_{k=1}^\infty \frac{\binom{4k}{k}}{k16^k}-214+\frac{40}{3}\log{2}.
	\end{align*}
	Applying \eqref{Fj,d}--\eqref{Fj,ddd} with $j=2$ and $x=\frac{1}{16}$, we get
	\begin{equation}
		A_2=\int_{1}^{\alpha}\frac{2(y-\alpha)}{y((3y+1)(y-1)-4y^2\sqrt{\frac{y-1}{3y+1}})}dy,
	\end{equation}
	\begin{equation}
		 A_2'=\int_{1}^{\alpha}\frac{2(4(3y^3-2y^2-y)-\al')}{y((3y+1)(y-1)-4y^2\sqrt{\frac{y-1}{3y+1}})}dy,
	\end{equation}
	\begin{equation}
		 A_2''=\int_{1}^{\alpha}\frac{2(16(27y^5-30y^4-16y^3+14y^2+5y)-\al'')}{y((3y+1)(y-1)-4y^2\sqrt{\frac{y-1}{3y+1}})}dy.
	\end{equation}
	So we deduce further that
	\begin{equation}\label{sigma2}
		\begin{aligned}
			\sigma_2&=I_2+\frac{287}{16}\frac{64\alpha^5}{(3\alpha+1)^2}-115(\alpha-1)
			-\frac{50}{3}\log\frac{4\alpha}{3\alpha+1}-214+\frac{40}{3}\log{2}
		\end{aligned}
	\end{equation}
	with the aid of the identity \eqref{pqr,1/16}, where
	 $$I_2:=\int_{1}^{\alpha}\frac{(27y^2-3y-40)(11y^3-11y^2-7y-1)}{4y((3y+1)(y-1)-4y^2\sqrt{\frac{y-1}{3y+1}})}dy.$$
	
	Putting $x=1/16$ in (\ref{equation of f}), we get $\sqrt{\frac{\alpha-1}{3\alpha+1}}=\frac{4\alpha^2}{(3\alpha+1)^2}$.
	For $z=\sqrt{(y-1)/{(3y+1)}}$, clearly
	$$y=-\f{z^2+1}{3z^2-1}\ \ \t{and}\ \ dy=\frac{8z}{(3z^2-1)^2}dz.$$
	So
	\begin{equation*}
		\begin{aligned}I_2
			 &=\int_{0}^{\sqrt{\frac{\alpha-1}{3\alpha+1}}}\frac{16(z^3-z^2+3z+1)(81z^4-75z^2+4)}{(z-1)(3z^2-1)^4 (z^2+1)}dz\\
			 &=\int_{0}^{\frac{4\alpha^2}{(3\alpha+1)^2}}\frac{16(z^3-z^2+3z+1)(81z^4-75z^2+4)}{(z-1)(3z^2-1)^4 (z^2+1)}dz.
		\end{aligned}
	\end{equation*}
	
	For the function
	\begin{align*}
		g_2(z)=\frac{4z(486z^5+135z^4-351z^3-126z^2+27z+11)}{(3z^2-1)^3}+10\log\frac{(z-1)^2}{z^2+1},
	\end{align*}
	it is easy to verify that
	$$g_2'(z)=\frac{16(z^3-z^2+3z+1)(81z^4-75z^2+4)}{(z-1)(3z^2-1)^4 (z^2+1)}.$$
	Thus,
	$$I_2=g_2\l(\frac{4\alpha^2}{(3\alpha+1)^2}\r)-g_2(0).$$
	Combining this with \eqref{sigma2}, we obtain
	\begin{align*}\sigma_2&=\frac{(11\alpha^3-11\alpha^2-7\alpha-1)p_1(\alpha)}{(3\alpha+1)^2 (33\alpha^4+108\alpha^3+54\alpha^2+12\alpha+1)^3}+\frac{10}{3}\log\frac{(\alpha+1)^6(3\alpha+1)^5(5\alpha+1)^6}{64\alpha^5(97\alpha^4+ 108 \alpha^3+ 54\alpha^2+12\alpha+1)^3}\\
		 &=\frac{10}{3}\log\l(1-\frac{(11\alpha^3-11\alpha^2-7\alpha-1)p_2(\alpha)}{64\alpha^5(97\alpha^4+ 108 \alpha^3+ 54\alpha^2+12\alpha+1)^3}\r)
		=0,
	\end{align*}
	where
	\begin{align*}
		p_1(y)&=3750516 y^{14}+40573764 y^{13}+ 178503183 y^{12}+ 437065524 y^{11}+ 654760746 y^{10}\\
		&\ \ \ +649292868 y^9+447600789 y^8+ 221967000 y^7+ 80835264 y^6+ 21786588 y^5   \\
		&\ \ \ +4319865 y^4+ 615636 y^3+ 59938 y^2+ 3580 y+99,
	\end{align*}
	and
	\begin{align*}
		p_2(y)&=4964927 y^{14}+ 19641236 y^{13}+ 39021903 y^{12}+ 50600432 y^{11}
		\\&\ \ \ + 47083971 y^{10}+ 32761508 y^9 +17347827 y^8+ 7024768 y^7+ 2167757 y^6
		\\&\ \ \ + 504220 y^5+ 86653 y^4+ 10640 y^3+ 881 y^2+ 44 y+1.
	\end{align*}
	This proves (\ref{H_2k}).
	\qed

	\medskip
	\noindent{\tt Proof of  $(\ref{H_3k})$}.
	Let
	$$\sigma_3:=\sum_{k=1}^\infty \frac{\binom{4k}{k}}{16^k}\l(P(k)H_{3k}-296k+178-\frac{25}{3k}\r)+196+\frac{80}{3}\log{2}.$$
	Taking the derivative and the second-order derivative of $F_3(x)$, we deduce that
	\begin{align*}
		\sigma_3 =\frac{22}{256}A_3''-\frac{70}{16}A_3'+11A_3-\frac{296}{16}\al'+178(\al-1)-\frac{25}{3}\sum_{k=1}^\infty \frac{\binom{4k}{k}}{k16^k}+196+\frac{80}{3}\log{2}.
	\end{align*}
	
	Applying \eqref{Fj,d}--\eqref{Fj,ddd} with $j=3$ and $x=1/{16}$, we get
	\begin{equation}
		A_3=\int_{1}^{\alpha}\frac{4(y-\al)}{3y(y-1)(3y+1-2\sqrt[3]{2} y\sqrt[3]{\frac{y}{y-1}})}dy, \label{F3,d}
	\end{equation}
	\begin{equation}	 A_3'=\int_{1}^{\alpha}\frac{4(4(3y^3-2y^2-y)-\al')}{3y(y-1)(3y+1-2\sqrt[3]{2} y\sqrt[3]{\frac{y}{y-1}})}dy,\label{F3,dd}
	\end{equation}
	\begin{equation}
		A_3''=\int_{1}^{\alpha}\frac{4(16(27y^5-30y^4-16y^3+14y^2+5y)-\al'')}{3y(y-1)(3y+1-2\sqrt[3]{2} y\sqrt[3]{\frac{y}{y-1}})}dy. \label{F3,ddd}
	\end{equation}
	So we deduce further that
	\begin{equation}\label{sigma3}
		\begin{aligned}
			\sigma_3&=I_3 -\frac{296}{16}\cdot\frac{64\alpha^5}{(3\alpha+1)^2}+178(\alpha-1)
			-\frac{100}{3}\log\frac{4\alpha}{3\alpha+1}+196+\frac{80}{3}\log{2}
		\end{aligned}
	\end{equation}
	with the aid of the identity \eqref{pqr,1/16}, where
	$$I_3:=\int_{1}^{\alpha}\frac{(27y^2-3y-40)(11y^3-11y^2-7y-1)}{6y(y-1)(3y+1-2\sqrt[3]{2} y\sqrt[3]{\frac{y}{y-1}})}dy.$$
	
	Putting $x=1/16$ in (\ref{equation of f}), we get $\sqrt[3]{\frac{\alpha-1}{\alpha}}=\frac{2\sqrt[3]{2}\alpha}{3\alpha+1}$.
	For $z=\sqrt[3]{(y-1)/{y}}$, clearly
	$$y=\f1{1-z^3}\ \ \t{and}\ \ dy=\frac{3z^2}{(z^3-1)^2}dz.$$
	So
	\begin{equation*}
		\begin{aligned}I_3
			 &=\int_{0}^{\sqrt[3]{\frac{\alpha-1}{\alpha}}}\frac{(40z^6-83z^3+16)(z^9-10z^6+28z^3-8)}{2(z^3-1)^4(z^4-4z+2\sqrt[3]{2})}dz\\
			 &=\int_{0}^{\frac{2\sqrt[3]{2}\alpha}{3\alpha+1}}\frac{(40z^6-83z^3+16)(z^9-10z^6+28z^3-8)}{2(z^3-1)^4(z^4-4z+2\sqrt[3]{2})}dz.
		\end{aligned}
	\end{equation*}
	
	For the function
	\begin{align*}
		 g_3(z)&=\frac{z(72z^8+40\sqrt[3]{2}z^7+20\sqrt[3]{4}z^6-135z^5-80\sqrt[3]{2}z^4-44\sqrt[3]{4}z^3+36z^2+22\sqrt[3]{2}z+12\sqrt[3]{4})}{2(z^3-1)^3}\\
		&\ \ \ -\frac{20}{3}\log\frac{2}{(\sqrt[3]{2}-z)^3},
	\end{align*}
	it is easy to verify that
	$$g_3'(z)=\frac{(40z^6-83z^3+16)(z^9-10z^6+28z^3-8)}{2(z^3-1)^4(z^4-4z+2\sqrt[3]{2})}.$$
	Thus,
	$$I_3=g_3\l(\frac{2\sqrt[3]{2}\alpha}{3\alpha+1}\r)-g_3(0).$$
	Combining this with \eqref{sigma3}, we obtain
	\begin{align*}\sigma_3&=-\frac{2(11\alpha^3-11\alpha^2-7\alpha-1)p_3(\alpha)}{(3\alpha+1)^2 (11\alpha^3+27\alpha^2+9\alpha+1)^3}
		\\&\ \ \ -\frac{20}{3}
		 \log\l(1+\frac{(11\alpha^3-11\alpha^2-7\alpha-1)(5\alpha^2+2\alpha+1)}{(\alpha+1)^3(3\alpha+1)^2}\r)\\
		&=0,
	\end{align*}
	where
	\begin{align*}
		p_3(y)&=71632y^{11}+ 599104y^{10}+ 2018295y^9+ 3436765y^8+3356404y^7 +2077344y^6\\
		&\ \ \ + 845834y^5+ 229586y^4+ 41220y^3+ 4712y^2+ 311y+9 .
	\end{align*}
	This proves (\ref{H_3k}).
	\qed
	
	\medskip
	\noindent{\tt Proof of  $(\ref{H_4k})$}.
	Let
	$$\sigma_4:=\sum_{k=1}^\infty \frac{\binom{4k}{k}}{16^k}\l(P(k)H_{4k}-\frac{449}{2}k+\frac{275}{2}-\frac{85}{12k}\r)+151+\frac{80}{3}\log{2}.$$
	Taking the derivative and the second-order derivative of $F_3(x)$, we deduce that
	\begin{align*}
		\sigma_4
		 =\frac{22}{256}A_4''-\frac{70}{16}A_4'+11A_4-\frac{449}{32}\al'+\frac{275}{2}(\al-1)-\frac{85}{12}\sum_{k=1}^\infty \frac{\binom{4k}{k}}{k16^k}+151+\frac{80}{3}\log{2}.
	\end{align*}
	
	Applying \eqref{Fj,d}--\eqref{Fj,ddd} with $j=4$ and $x=\frac{1}{16}$, we get
	\begin{equation}
		A_4=\int_{1}^{\alpha}\frac{y-\alpha}{y(y-1)(3y+1-2y\sqrt[4]{\frac{3y+1}{y-1}})}dy,\label{F4,d}
	\end{equation}
	\begin{equation}
		 A_4'=\int_{1}^{\alpha}\frac{4(3y^3-2y^2-y)-\al'}{y(y-1)(3y+1-2y\sqrt[4]{\frac{3y+1}{y-1}})}dy,\label{F4,dd}
	\end{equation}
	\begin{equation}
		 A_4''=\int_{1}^{\alpha}\frac{16(27y^5-30y^4-16y^3+14y^2+5y)-\al''}{y(y-1)(3y+1-2y\sqrt[4]{\frac{3y+1}{y-1}})}dy.\label{F4,ddd}
	\end{equation}
	
	So we deduce further that
	\begin{equation}\label{sigma4}
		\begin{aligned}
			\sigma_4&=I_4-\frac{449}{32}\cdot\frac{64\alpha^5}{(3\alpha+1)^2}+\frac{275}{2}(\alpha-1)
			-\frac{85}{3}\log\frac{4\alpha}{3\alpha+1}+151+\frac{80}{3}\log{2}
		\end{aligned}
	\end{equation}
	with the aid of the identity \eqref{pqr,1/16}, where
	 $$I_4:=\int_{1}^{\alpha}\frac{(27y^2-3y-40)(11y^3-11y^2-7y-1)}{8y(y-1)(3y+1-2y\sqrt[4]{\frac{3y+1}{y-1}})}dy
	.$$
	
	Putting $x=1/16$ in (\ref{equation of f}), we get $\sqrt[4]{\frac{\alpha-1}{3\alpha+1}}=\frac{2\alpha}{3\alpha+1}$.
	For $z=\sqrt[4]{(y-1)/{(3y+1)}}$, clearly
	$$y=-\f{z^4+1}{3z^4-1}\ \ \t{and}\ \ dy=\frac{16z^3}{(3z^4-1)^2}dz.$$
	So
	\begin{equation*}
		\begin{aligned}I_4
			&=\int_{0}^{\sqrt[4]{\frac{\alpha-1}{3\alpha+1}}}\frac{8(z^3-z^2+z+1)(z^6-z^4+3z^2+1) (81z^8-75z^4+4)}{(z-1)(3z^4-1)^4(z^4+1)}dz\\
			 &=\int_{0}^{\frac{2\alpha}{3\alpha+1}}\frac{8(z^3-z^2+z+1)(z^6-z^4+3z^2+1)(81z^8-75z^4+4)}{(z-1) (3z^4-1)^4(z^4+1)}dz.
		\end{aligned}
	\end{equation*}
	
	For the function
	\begin{align*} g_4(z)&=\frac{2zQ(z)}{(3 z^4-1)^3}+5\log\frac{(z-1)^4}{z^4+1}
	\end{align*}
	with
	$$Q(z)=486z^{11}+270z^{10}+135z^9+54z^8-351z^7-216z^6-126z^5-68z^4+27z^3+18z^2+11z+6,$$
	it is easy to verify that
	$$g_4'(z)=\frac{8(z^3-z^2+z+1)(z^6-z^4+3z^2+1) (81z^8-75z^4+4)}{(z-1) (3z^4-1)^4(z^4+1)}.$$
	Thus,
	$$I_4=g_4\l(\frac{2\alpha}{3\alpha+1}\r)-g_4(0).$$
	Combining this with \eqref{sigma4} and
	the identity $$\frac{(3\alpha+1)^{15}}{\alpha^{15}}=\l(\frac{16\alpha}{\alpha-1}\r)^5$$
	obtained from (\ref{equation of f}) with $x=1/{16}$, we find that
	\begin{align*}\sigma_4&=-\frac{(11\alpha^3-11\alpha^2-7\alpha-1)p_4(\alpha)}{2(1+3\alpha)^2 (11\alpha^3+27\alpha^2+9\alpha+1)^3}
		\\&\ \ \ +\frac{5}{3}\log\frac{(\alpha+1)^{12}(3\alpha+1)^{17}}{262144\alpha^{17}(97\alpha^4+108\alpha^3+54\alpha^2+ 12\alpha+1)^3}\\
		&=\frac{5}{3}\log\frac{16^5\alpha^3(\alpha+1)^{12} (3\alpha+1)^2}{262144(\alpha-1)^5 (97\alpha^4+108\alpha^3+54\alpha^2+ 12\alpha+1)^3}\\
		&=\frac{5}{3}\log\l(1-\frac{(11\alpha^3-11\alpha^2-7\alpha-1)p_5(\alpha)}{(\alpha-1)^5 (97\alpha^4+108\alpha^3+54\alpha^2+ 12\alpha+1)^3}\r)
		=0,
	\end{align*}
	where
	\begin{align*} p_4(y)&=5867532y^{14}+63476028y^{13}+276465231y^{12}+604283868y^{11}+773801154y^{10}
		\\&\ \ \ +655999884y^9+390925197y^8+169175304y^7+54045864y^6+12792132y^5
		\\&\ \ \ +2221641y^4+275580y^3+23114y^2+1172y+27
	\end{align*}
	and
	\begin{align*}
		p_5(y)&=82967 y^{14}- 54788 y^{13}- 111085y^{12}+ 20604 y^{11}+ 108295 y^{10}+ 41488 y^9\\
		&\ \ \ -30309y^8- 30728y^7 - 5139y^6 +5660y^5+ 4177y^4+ 1356y^3+ 245y^2+ 24y+1 .	
	\end{align*}
	This proves (\ref{H_4k}).
	\qed
	
	\section{Proof of Theorem $\ref{2}$}
	\setcounter{lemma}{0}
	\setcounter{theorem}{0}
	\setcounter{equation}{0}
	\setcounter{conjecture}{0}
	\setcounter{remark}{0}
	\setcounter{corollary}{0}
	
	\medskip
	
	\begin{lemma}\label{lemma4.1}
		Let $m$ be any nonzero complex number and let $\psi$ be a  function from $\N$ to $\C$.
		Let $\Delta\psi(k)=\psi(k+1)-\psi(k)$ for $k\in\Z^+$. Then
		\begin{equation}
			\begin{aligned}
				&\sum_{k=1}^{n}\frac{\binom{4k}{k}\psi(k)((256-27m)k^3+384k^2+(176+3m)k+24 )}{(3k+1)m^k}\\
				=&\frac{8(2n+1)(4n+1)(4n+3)\binom{4n}{n}\psi(n)}{(3n+1)m^n}
				-\sum_{k=0}^{n-1}\frac{8(2k+1)(4k+1)(4k+3)\binom{4k}{k}}{(3k+1)m^k}\Delta \psi(k)		 
			\end{aligned}\label{H_jk,3k+1}
		\end{equation}
		and
		\begin{equation}
			\begin{aligned}
				&\ \sum_{k=1}^{n}\frac{\binom{4k}{k}\psi(k)((256-27m)k^3+3(128-9m)k^2+2(88-3m)k+24 )}{(3k+1)(3k+2)m^k}\\
				=&\ \frac{8(2n+1)(4n+1)(4n+3)\binom{4n}{n}\psi(n)}{(3n+1)(3n+2)m^n}
				 -\sum_{k=0}^{n-1}\frac{8(2k+1)(4k+1)(4k+3)\binom{4k}{k}}{(3k+1)(3k+2)m^k}\Delta\psi(k).
			\end{aligned}\label{H_jk,3k+1,3k+2}
		\end{equation}
	\end{lemma}
	\Proof. This can be easily proved by induction or the Abel summation. \qed
	
	\begin{lemma} For any $m,n\in\Z^+$, we have
		\begin{equation}
			G_m(x)^n=1+\sum_{k=1}^{\infty}\frac{n}{k}\binom{mk+n-1}{k-1}x^k.
			\label{Gm(x)^s}
		\end{equation}
	\end{lemma}
	\Proof. Let $A(x)=G_m(x)-1,\Phi(x)=(1+x)^m$. Then $A(x)=x\Phi(A(x))$.
	Applying Lagrange's inversion formula (cf.\cite{7}), for each $k=1,2,3,\ldots$ we have
	\begin{equation}
		k[x^k]G_m(x)^n=k[x^k]H(A(x))=[z^{k-1}]H'(z)\Phi(z)^k=[z^{k-1}]n(1+z)^{mk+n-1},
	\end{equation}
	where $H(z)=(1+z)^n$, and $[x^k]\Psi(x)$ denotes the coefficient of $x^k$
	in the power series expansion of $\Psi(x)$.
	This proves \eqref{Gm(x)^s}. \qed
	
	\begin{remark} Applying \eqref{Gm(x)^s} with $m=4$ and $n\in\{2,3\}$, we immediately get the identities
		\begin{equation}
			 \sum_{k=0}^{\infty}\frac{\binom{4k}{k}(8k+2)x^k}{(3k+1)(3k+2)}=\l(\frac{4f(x)}{3f(x)+1}\r)^2
			\label{f,5}
		\end{equation}
		and
		\begin{equation}
			 \sum_{k=0}^\infty\frac{\binom{4k}{k}(4k+1)(4k+2)}{(k+1)(3k+1)(3k+2)}x^k=\l(\frac{4f(x)}{3f(x)+1}\r)^3.
			\label{f,6}
		\end{equation}
	\end{remark}

	\begin{lemma}
		We have
		\begin{equation} \sum_{k=1}^\infty\frac{\binom{4k}{k}(638k^4+966k^3+194k^2-199k-59)}{(k+1)(3k+1)(3k+2)16^k}=\frac{103}{2}
			\label{103}
		\end{equation}
		and
		\begin{equation}	 \sum_{k=1}^\infty\frac{\binom{4k}{k}(770k^4+1134k^3+82k^2-381k-105)}{(k+1)(3k+1)(3k+2)16^k}=\frac{117}{2}.
			\label{117}
		\end{equation}
	\end{lemma}
	\Proof.
	Let $\alpha=f(\frac{1}{16})$.
	Putting $x=1/16$ in (\ref{f,1}), (\ref{f,5}) and (\ref{f,6}), we find
	\begin{align}
		\sum_{k=0}^\infty\frac{\binom{4k}{k}}{(3k+1)16^k}&=\frac{4\alpha}{3\alpha+1},
		\\
		 \sum_{k=0}^\infty\frac{\binom{4k}{k}(8k+2)}{(3k+1)(3k+2)16^k}&=\l(\frac{4\alpha}{3\alpha+1}\r)^2,
		\\ \sum_{k=0}^\infty\frac{\binom{4k}{k}(4k+1)(4k+2)}{(k+1)(3k+1)(3k+2)16^k}&
		=\l(\frac{4\alpha}{3\alpha+1}\r)^3.
	\end{align}
	Thus, with the aid of \eqref{f,3} we deduce that
	\begin{align*}
		&\ \ \ \sum_{k=1}^\infty\frac{\binom{4k}{k}(638k^4+966k^3+194k^2-199k-59)}{(k+1)(3k+1)(3k+2)16^k}\\
		 &=\sum_{k=1}^\infty\frac{\binom{4k}{k}(4k+1)(4k+2)}{(k+1)(3k+1)(3k+2)16^k}+\frac{23}{18}\sum_{k=1}^\infty\frac{\binom{4k}{k}(8k+2)}{(3k+1)(3k+2)16^k} \\
		&\ \ \ +\frac{8}{3}\sum_{k=1}^\infty\frac{\binom{4k}{k}}{(3k+1)16^k}
		+\frac{638}{9}\sum_{k=1}^\infty\frac{\binom{4k}{k}k}{16^k}
		-\frac{310}{9}\sum_{k=1}^\infty\frac{\binom{4k}{k}}{16^k}\\
		&=\l(\frac{4\alpha}{3\alpha+1}\r)^3-1+\frac{23}{18}\l(\l(\frac{4\alpha}{3\alpha+1}\r)^2-1\r)
		+\frac{8}{3}\l(\frac{4\alpha}{3\alpha+1}-1\r)
		\\&\ \ \ +\frac{638}{9\times16}\cdot\frac{64\alpha^5}
		{(3\alpha+1)^2}-\frac{310}{9}(\alpha-1)\\
		&=\frac{2(11\alpha^3-11\alpha^2- 7\alpha-1)(348 \alpha^3+ 464\alpha^2+ 305\alpha+99)}{9(3 \alpha+1)^3}+\frac{103}{2}=\f{103}2.
	\end{align*}
	Similarly, we have
	\begin{align*}
		&\ \ \ \sum_{k=1}^\infty\frac{\binom{4k}{k}(770k^4+1134k^3+82k^2-381k-105)}{(k+1)(3k+1)(3k+2)16^k}\\
		&=-\sum_{k=1}^\infty\frac{\binom{4k}{k}(4k+1)(4k+2)}{(k+1)(3k+1)(3k+2)16^k}
		-\frac{43}{18}\sum_{k=1}^\infty\frac{\binom{4k}{k}(8k+2)}{(3k+1)(3k+2)16^k} \\
		&\ \ \ -4\sum_{k=1}^\infty\frac{\binom{4k}{k}}{(3k+1)16^k}+\frac{770}{9}
		\sum_{k=1}^\infty\frac{\binom{4k}{k}k}{16^k}
		-\frac{406}{9}\sum_{k=1}^\infty\frac{\binom{4k}{k}}{16^k}\\
		&=-\l(\frac{4\alpha}{3\alpha+1}\r)^3+1-\frac{43}{18}\l(\l(\frac{4\alpha}{3\alpha+1}\r)^2-1\r)
		-4\l(\frac{4\alpha}{3\alpha+1}-1\r)
		\\&\ \ \ +\frac{770}{9\times16}\cdot\frac{64\alpha^5}
		{(3\alpha+1)^2}-\frac{406}{9}(\alpha-1)\\
		&=\frac{2(11\alpha^3-11\alpha^2- 7\alpha-1)(420 \alpha^3+ 560 \alpha^2+ 329 \alpha+27)}{9(3 \alpha+1)^3}+\frac{117}{2}=\f{117}2.
	\end{align*}
	This concludes the proof. \qed

	\begin{lemma}
		We have
		\begin{equation}
			\sum_{k=1}^\infty \frac{\binom{4k}{k}((11k^2+8k+1)H_{2k}+\frac{23}{2}k+\frac{17}{2}+\frac{5}{3k})}{(3k+1)(3k+2)16^k}
			=\frac{8}{3}\log2-\frac{1}{2}.
			\label{H_2k,16,(3k+1)(3k+2)}
		\end{equation}
	\end{lemma}
	\Proof. Observe that
	\begin{equation}\begin{aligned}
			 \frac{11k^2+8k+1}{(3k+1)(3k+2)}&=-\frac{22k^2-92k+11}{5}-\frac{3(-176k^3-48k^2+80k+24)}{40(3k+1)}
			\\&\ \ \ +\frac{3(-176k^3+384k^2+224k+24)}{8(3k+1)(3k+2)}.
		\end{aligned}\label{kk}\end{equation}
	
	Putting $m=16$ and $\psi(k)=H_{2k}$ in (\ref{H_jk,3k+1}) and (\ref{H_jk,3k+1,3k+2}) and letting $n\to +\infty$, we then obtain $S_1=0=S_2$, where $S_1$ denotes the expression
	\begin{equation*} \sum_{k=1}^{\infty}\frac{\binom{4k}{k}H_{2k}(-176k^3+384k^2+224k+24)}{(3k+1)16^k}
		+\sum_{k=0}^{\infty}\frac{8(2k+1)(4k+1)(4k+3)\binom{4k}{k}}{(3k+1)16^k}(H_{2k+2}-H_{2k})
	\end{equation*}
	and $S_2$ denotes the expression
	\begin{equation*} \sum_{k=1}^{\infty}\frac{\binom{4k}{k}H_{2k}(-176k^3-48k^2+80k+24)}{(3k+1)(3k+2)16^k}
		 +\sum_{k=0}^{\infty}\frac{8(2k+1)(4k+1)(4k+3)\binom{4k}{k}}{(3k+1)(3k+2)16^k}(H_{2k+2}-H_{2k}).
	\end{equation*}
	
	In light of \eqref{kk} and the fact that $S_1=0=S_2$, we have
	\begin{align*}
		&\ \ \ \frac{1}{5}\sum_{k=1}^\infty\frac{\binom{4k}{k}}{16^k}\l(P(k)H_{2k}+287k-115-\frac{25}{6k}\r)
		+\sum_{k=1}^\infty \frac{\binom{4k}{k}((11k^2+8k+1)H_{2k}+\frac{23}{2}k+\frac{17}{2}+\frac{5}{3k})}{(3k+1)(3k+2)16^k}
		\\
		&=\frac{1}{5}\sum_{k=1}^\infty\frac{\binom{4k}{k}(287k-115-\frac{25}{6k})}{16^k}
		+\frac{3}{40}\l(\sum_{k=1}^\infty\frac{8(2k+1)(4k+1)(4k+3)\binom{4k}{k}}{(3k+1)16^k}
		(H_{2k+2}-H_{2k})+36\r)\\
		&\ \ \ -\frac{3}{8}\l(\sum_{k=1}^\infty\frac{8(2k+1)(4k+1)(4k+3)\binom{4k}{k}}{(3k+1)(3k+2)16^k}
		 (H_{2k+2}-H_{2k})+18\r)+\sum_{k=1}^\infty\frac{\binom{4k}{k}\l(\frac{23}{2}k+\frac{17}{2}+\frac{5}{3k}\r)}{(3k+1)(3k+2)16^k}
		\\
		 &=\frac{9}{10}\sum_{k=1}^\infty\frac{\binom{4k}{k}(638k^4+966k^3+194k^2-199k-59)}{(k+1)(3k+1)(3k+2)16^k}
		-\frac{81}{20}=\f{423}{10}
	\end{align*}
	with the aid of \eqref{103}.
	Thus, we have
	\begin{align*}
		&\ \sum_{k=1}^\infty \frac{\binom{4k}{k}((11k^2+8k+1)H_{2k}+\frac{23}{2}k+\frac{17}{2}+\frac{5}{3k})}{(3k+1)(3k+2)16^k}
		\\=&\ -\frac{1}{5}\sum_{k=1}^\infty\frac{\binom{4k}{k}}{16^k}\l(P(k)H_{2k}+287k-115-\frac{25}{6k}\r)+\frac{423}{10}.
	\end{align*}
	So \eqref{H_2k,16,(3k+1)(3k+2)} follows from \eqref{H_2k}.\qed
	
	\begin{lemma}
		We have
		\begin{equation}
			 \sum_{k=1}^\infty\frac{\binom{4k}{k}\l((11k^2+8k+1)H_{4k}+\frac{79}{4}k+\frac{65}{4}+\frac{17}{6k}\r)}{(3k+1)(3k+2)16^k}
			=\frac{16}{3}\log2-\frac{7}{4}.
			\label{H_4k,16,(3k+1)(3k+2)}
		\end{equation}
	\end{lemma}
	
	\textit{Proof.} Putting $m=16$ and $\psi(k)=H_{4k}$ in \eqref{H_jk,3k+1} and \eqref{H_jk,3k+1,3k+2}, and letting $n\to +\infty$, we then obtain $S_3=0=S_4$, where $S_3$ denotes the expression
	\begin{equation*}	 \sum_{k=1}^{\infty}\frac{\binom{4k}{k}H_{4k}(-176k^3+384k^2+224k+24)}{(3k+1)16^k}
		 +\sum_{k=0}^{\infty}\frac{8(2k+1)(4k+1)(4k+3)\binom{4k}{k}(H_{4k+4}-H_{4k})}{(3k+1)16^k}\label{16,H_4k,3k+1}
	\end{equation*}
	and $S_4$ stands for
	\begin{equation*}	 \sum_{k=1}^{\infty}\frac{\binom{4k}{k}H_{4k}(-176k^3-48k^2+80k+24)}{(3k+1)(3k+2)16^k}+\sum_{k=0}^{\infty}\frac{8(2k+1)(4k+1)(4k+3)\binom{4k}{k}(H_{4k+4}-H_{4k})}{(3k+1)(3k+2)16^k}.\label{16,H_4k,3k+1,3k+2}
	\end{equation*}
	
	In light of \eqref{kk} and the facts that $S_3=0=S_4$, we have
	\begin{align*}
		&\ \ \ \frac{1}{5}\sum_{k=1}^\infty \frac{\binom{4k}{k}}{16^k}\l(P(k)H_{4k}-\frac{449k-275}{2}-\frac{85}{12k}\r)
		 +\sum_{k=1}^\infty\frac{\binom{4k}{k}\big((11k^2+8k+1)H_{4k}+\frac{79}{4}k+\frac{65}{4}+\frac{17}{6k}\big)}{(3k+1)(3k+2)16^k}\\
		&=\frac{1}{5}\sum_{k=1}^\infty\frac{\binom{4k}{k}(-\frac{449k-275}{2}-\frac{85}{12k})}{16^k} +\frac{3}{40}\l(\sum_{k=1}^\infty\frac{8(2k+1)(4k+1)(4k+3)\binom{4k}{k}}{(3k+1)16^k}(H_{4k+4}-H_{4k})+50\r)\\
		&\ \ \ +\sum_{k=1}^\infty\frac{\binom{4k}{k}\l(\frac{79}{4}k+\frac{65}{4}+\frac{17}{6k}\r)}{(3k+1)(3k+2)16^k}-\frac{3}{8}\l(\sum_{k=1}^\infty\frac{8(2k+1)(4k+1)(4k+3)\binom{4k}{k}}{(3k+1)(3k+2)16^k}(H_{4k+4}-H_{4k})+25\r)\\
		 &=-\frac{9}{20}\sum_{k=1}^\infty\frac{\binom{4k}{k}(770k^4+1134k^3+82k^2-381k-105)}{(1+k)(1+3k)(2+3k)16^k}-\frac{45}{8}\\
		&=-\frac{9}{20}\times\frac{117}{2}-\frac{45}{8}=-\frac{639}{20}
	\end{align*}
	with the aid of \eqref{117}.
	Thus, we have
	\begin{align*}
		&\ \sum_{k=1}^\infty \frac{\binom{4k}{k}((11k^2+8k+1)H_{4k}+\frac{79}{4}k+\frac{65}{4}+\frac{17}{6k})}{(3k+1)(3k+2)16^k}
		\\=&\ -\frac{1}{5}\sum_{k=1}^\infty \frac{\binom{4k}{k}}{16^k}\l(P(k)H_{4k}-\frac{449k-275}{2}-\frac{85}{12k}\r)-\frac{639}{20}.
	\end{align*}
	So \eqref{H_4k,16,(3k+1)(3k+2)} follows from \eqref{H_4k}.\qed
	
	\medskip
	\noindent{\bf Proof of Theorem \ref{2}}.
	Putting $m=16$ and $\psi(k)=H_k$ in \eqref{H_jk,3k+1} and \eqref{H_jk,3k+1,3k+2} and
	letting $n\to +\infty$, we then obtain
	\begin{equation}
		 \sum_{k=1}^{\infty}\frac{\binom{4k}{k}\big(H_k(-176k^3+384k^2+224k+24)+48(3k-1)(3k+1)\big)}{(3k+1)16^k}
		=0\label{16,H_k,3k+1}
	\end{equation}
	and
	\begin{equation}
		\sum_{k=1}^{\infty}\frac{\binom{4k}{k}\big(H_k(-176k^3-48k^2+80k+24)+48(3k+1)(3k+2)\big)}
		{(3k+1)(3k+2)16^k}=0.\label{16,H_k,3k+1,3k+2}
	\end{equation}
	
	In light of \eqref{kk}, \eqref{16,H_k,3k+1} and \eqref{16,H_k,3k+1,3k+2}, we have
	\begin{align*}
		&\ \ \ \frac{1}{5}\sum_{k=1}^\infty \frac{\binom{4k}{k}}{16^k}\l(P(k)H_{k}-54k+108-\frac{10}{3k}\r)+\sum_{k=1}^\infty\frac{\binom{4k}{k}\l((11k^2+8k+1)H_k+6k+6+\frac{4}{3k}\r)}{(3k+1)(3k+2)16^k}\\
		&=\frac{1}{5}\sum_{k=1}^\infty\frac{\binom{4k}{k}(-54k+108-\frac{10}{3k})}{16^k}
		 +\frac{3}{40}\sum_{k=1}^\infty\frac{48(3k-1)\binom{4k}{k}}{16^k}-\frac{3}{8}\sum_{k=1}^\infty\frac{48\binom{4k}{k}}{16^k}\\
		&\ \ \ +\sum_{k=1}^\infty\frac{\bi{4k}k(6k+6+\frac{4}{3k})}{(3k+1)(3k+2)16^k}\\
		&=\sum_{k=1}^\infty \f{\binom{4k}{k}}{16^k}\l(\frac{1}{5}(-54k+108-\frac{10}{3k})+\frac{18}{5}(3k-1)-18-\frac{2}{k}\r)=0.
	\end{align*}
	So \eqref{H_k,16,(3k+1)(3k+2)} follows from \eqref{H_k}.
	
	Similarly, in light of \eqref{H_k,16,(3k+1)(3k+2)} and \eqref{H_2k,16,(3k+1)(3k+2)},  we have
	\begin{align*}
		&\ \ \ \sum_{k=0}^\infty\frac{\binom{4k}{k}((11k^2+8k+1)(H_{2k}-\frac{5}{4}H_k)+4k+1)}{(3k+1)(3k+2)16^k}-\frac{1}{2}\\
		 &=\sum_{k=1}^\infty\frac{\binom{4k}{k}\big((11k^2+8k+1)H_{2k}+\frac{23}{2}k+\frac{17}{2}+\frac{5}{3k}\big)}{(3k+1)(3k+2)16^k}
		 -\sum_{k=1}^\infty\frac{\binom{4k}{k}\frac{5}{4}((11k^2+8k+1)H_k+6k+6+\frac{4}{3k})}{(3k+1)(3k+2)16^k}\\
		&=\l(\frac{8}{3}\log2-\frac{1}{2}\r)-\frac{5}{4}\times \frac{4}{3}\log2=\log2-\frac{1}{2}
	\end{align*}
	Similarly, in light of \eqref{H_2k,16,(3k+1)(3k+2)} and \eqref{H_4k,16,(3k+1)(3k+2)},  we have
	\begin{align*}
		&\ \ \ \sum_{k=0}^\infty\frac{\binom{4k}{k}\big((11k^2+8k+1)(10H_{4k}-17H_{2k})+2k+18\big)}{(3k+1)(3k+2)16^k}\\
		 &=10\sum_{k=1}^\infty\frac{\binom{4k}{k}\big((11k^2+8k+1)H_{4k}+\frac{79}{4}k+\frac{65}{4}+\frac{17}{6k}\big)}{(3k+1)(3k+2)16^k} \\
		&\ \ \ -17\sum_{k=1}^\infty\frac{\binom{4k}{k}\big((11k^2+8k+1)H_{2k}+\frac{23}{2}k+\frac{17}{2}+\frac{5}{3k}\big)}{(3k+1)(3k+2)16^k}+9 \\
		&=10\times\l(\frac{16}{3}\log2-\frac{7}{4}\r)-17\times \l(\frac{8}{3}\log2-\frac{1}{2}\r)+9=8\log2.
	\end{align*}
	So \eqref{4H_2k-5H_k,16,(3k+1)(3k+2)} follows from \eqref{H_k}, and \eqref{H_2k,16,(3k+1)(3k+2)}. \eqref{10H_4k-17H_2k,16,(3k+1)(3k+2)} follows from \eqref{H_2k,16,(3k+1)(3k+2)} and \eqref{H_4k,16,(3k+1)(3k+2)}.
	
	In view of the above, we have completed the proof of Theorem \ref{2}.
	\qed

	\section{Proof of Theorem $\ref{3}$}
	\setcounter{lemma}{0}
	\setcounter{theorem}{0}
	\setcounter{equation}{0}
	\setcounter{conjecture}{0}
	\setcounter{remark}{0}
	\setcounter{corollary}{0}

	\begin{lemma}
		We have
		\begin{equation}	\begin{aligned}
				&\frac{64}{5}\sum_{k=0}^{\infty}\frac{\binom{4k}{k}(182k+5)}{(-256)^k}
				+\frac{36}{5}\sum_{k=0}^{\infty}\frac{\binom{4k}{k}}{(3k+1)(-256)^k}
				\\=&\ \sum_{k=0}^{\infty}\frac{\binom{4k}{k}(8k+2)}{(3k+1)(3k+2)(-256)^k}+\frac{72}{5}\sqrt{2},
			\end{aligned} \label{-256,lemma,equivalent}
		\end{equation}
		\begin{equation}\begin{aligned}
				&\ \frac{32}{5}\sum_{k=0}^{\infty}\frac{\binom{4k}{k}(725k-49)}{128^k}
				-\frac{12}{5}\sum_{k=0}^{\infty}\frac{\binom{4k}{k}}{(3k+1)128^k}
				\\=&\ -\sum_{k=0}^{\infty}\frac{\binom{4k}{k}(8k+2)}{(3k+1)(3k+2)128^k}-\frac{576}{5}\sqrt{2},
			\end{aligned}		\label{128,lemma,equivalent}
		\end{equation}
		\begin{equation}
			\begin{aligned}
				&\ \frac{1}{65}\sum_{k=0}^{\infty}\frac{\binom{4k}{k}242(175k+12)}{(-72)^k}
				-\frac{28}{5}\sum_{k=0}^{\infty}\frac{\binom{4k}{k}}{(3k+1)(-72)^k}
				\\=&\ \sum_{k=0}^{\infty}\frac{\binom{4k}{k}(8k+2)}{(3k+1)(3k+2)(-72)^k}+\frac{216}{65}\sqrt{3},
			\end{aligned}		\label{-72,lemma,equivalent}
		\end{equation}
		\begin{equation}
			\begin{aligned}
				&\ \frac{1}{460}\sum_{k=0}^{\infty}\frac{\binom{4k}{k}4(118237k+17320)}{(-25)^k}
				-\frac{96}{5}\sum_{k=0}^{\infty}\frac{\binom{4k}{k}}{(3k+1)(-25)^k}
				\\&=\ 4\sum_{k=0}^{\infty}\frac{\binom{4k}{k}(8k+2)}{(3k+1)(3k+2)(-25)^k}+\frac{72}{23}\sqrt{5},
			\end{aligned}	\label{-25,lemma,equivalent}
		\end{equation}
		\begin{equation}\begin{aligned}
				&\ \frac{1}{5}\sum_{k=0}^{\infty}\frac{\binom{4k}{k}(-3038k+1160)}{24^k}
				-\frac{4}{5}\sum_{k=0}^{\infty}\frac{\binom{4k}{k}}{(3k+1)24^k}
				 \\&=-\sum_{k=0}^{\infty}\frac{\binom{4k}{k}(8k+2)}{(3k+1)(3k+2)24^k}+\frac{216}{5}\sqrt{3}.
			\end{aligned}		\label{24,lemma,equivalent}
		\end{equation}
	\end{lemma}
	\Proof. We just prove \eqref{-256,lemma,equivalent} in details. Other formulas in the lemma can be proved similarly.
	
	Let $\beta=f(-1/{256})$. Applying \eqref{equation of f} with $x=-1/256$, we find that $$0=7(28\beta^4-18\beta^2-8\beta-1)=(14\beta^2-7\sqrt{2}\beta-1-2\sqrt{2})(14\beta^2+7\sqrt{2}
	\beta-1+2\sqrt{2}).$$  As $\beta\in \mathbb{R}$, we have $14\beta^2+7\sqrt{2}
	\beta-1+2\sqrt{2}\not=0$ and hence
	\begin{equation}
		14\beta^2-7\sqrt{2}\beta-1-2\sqrt{2}=0. \label{-256,eq}
	\end{equation}
	Thus, with the aids of \eqref{f,3}, \eqref{f,1}, \eqref{f,5} and \eqref{-256,eq}, we deduce that
	\begin{align*}
		&\ \ \ \frac{64}{5}\sum_{k=0}^{\infty}\frac{\binom{4k}{k}(182k+5)}{(-256)^k}-\frac{36}{5}\sum_{k=0}^{\infty}\frac{\binom{4k}{k}}{(3k+1)(-256)^k}
		-\sum_{k=0}^{\infty}\frac{\binom{4k}{k}(8k+2)}{(3k+1)(3k+2)(-256)^k}\\	 &=\frac{64}{5}\l(-\frac{182}{256}\cdot\frac{64\beta^5}{(3\beta+1)^2}+5\beta\r)-\frac{36}{5}\frac{4\beta}{3\beta+1}-\l(\frac{4\beta}{3\beta+1}\r)^2-\frac{72}{5}\sqrt{2}
		+\frac{72}{5}\sqrt{2}\\
		&=-\frac{8(14\beta^2-7\sqrt{2}\beta-1-2\sqrt{2})}{35(3\beta+1)^2}\lambda+\frac{72}{5}\sqrt{2}
		=\frac{72}{5}\sqrt{2},
	\end{align*}
	where
	\begin{equation*}
		\lambda=182\beta^3+ 91 \sqrt{2}\beta^2+ (-76 + 26 \sqrt{2})\beta-36 + 9 \sqrt{2}.
	\end{equation*}
	This proves the desired \eqref{-256,lemma,equivalent}. \qed
	
	\medskip
	\noindent{\tt Sketch of the Proof of Theorem \ref{3}}.
	Putting $m=-256$ and $\psi(k)=H(k)$ in \eqref{H_jk,3k+1} and \eqref{H_jk,3k+1,3k+2} and
	letting $n\to +\infty$, we then obtain
	\begin{equation}
		\sum_{k=0}^{\infty}\frac{\binom{4k}{k}(H(k)(896k^3+48k^2-74k+3)+2)}{(3k+1)(-256)^k}=0 \label{-256,1}
	\end{equation}
	and
	\begin{equation}	 \sum_{k=0}^{\infty}\frac{\binom{4k}{k}(H(k)(896k^3+912k^2+214k+3)+2)}{(3k+1)(3k+2)(-256)^k}=0.
		\label{-256,2}
	\end{equation}
	Observe that
	\begin{equation*}
		\begin{aligned}
			&\ \ \ \sum_{k=0}^\infty\frac{\binom{4k}{k}(112k^2+110k+23)H(k)}{(3k+1)(3k+2)(-256)^k}+\frac{64}{5}\sum_{k=0}^\infty\frac{\binom{4k}{k}(224k^2-86k+1)H(k)}{(-256)^k}\\
			 &=\frac{48}{5}\sum_{k=0}^{\infty}\frac{\binom{4k}{k}(896k^3+48k^2-74k+3)H(k)}{(3k+1)(-256)^k} -3\sum_{k=0}^{\infty}\frac{\binom{4k}{k}(896k^3+912k^2+214k+3)H(k)}{(3k+1)(3k+2)(-256)^k}.
		\end{aligned}
	\end{equation*}
	Combining this with \eqref{-256,1} and \eqref{-256,2}, we obtain
	\begin{align*}
		&\ \ \ \sum_{k=0}^\infty\frac{\binom{4k}{k}((112k^2+110k+23)H(k)+28k+16)}{(3k+1)(3k+2)(-256)^k}+
		\frac{64}{5}\sum_{k=0}^\infty\frac{\binom{4k}{k}(224k^2-86k+1)H(k)}{(-256)^k}\\
		&=\sum_{k=0}^\infty \frac{\binom{4k}{k}(28k+16)}{(3k+1)(3k+2)(-256)^k}-\frac{96}{5}\sum_{k=0}^{\infty}\frac{\binom{4k}{k}}{(3k+1)(-256)^k}+6\sum_{k=0}^{\infty}\frac{\binom{4k}{k}}{(3k+1)(3k+2)(-256)^k}\\
		 &=-\frac{36}{5}\sum_{k=0}^{\infty}\frac{\binom{4k}{k}}{(3k+1)(-256)^k}-\sum_{k=0}^{\infty}\frac{\binom{4k}{k}(8k+2)}{(3k+1)(3k+2)(-256)^k}\\
		&=\frac{72}{5}\sqrt{2}-\frac{64}{5}\sum_{k=0}^{\infty}\frac{\binom{4k}{k}(182k+5)}{(-256)^k}.
	\end{align*}
	with the aid of \eqref{-256,lemma,equivalent}.
	Thus
	\begin{align*}
		&\ \ \ \sum_{k=0}^\infty\frac{\binom{4k}{k}((112k^2+110k+23)H(k)+28k+16)}{(3k+1)(3k+2)(-256)^k}\\
		 &=\frac{72}{5}\sqrt{2}-\frac{64}{5}\sum_{k=0}^\infty\frac{\binom{4k}{k}((224k^2-86k+1)H(k)+182k+5)}{(-256)^k}.\label{-256,lemma1,h}
	\end{align*}
	and hence \eqref{-256,3k+2,h} is equivalent to \eqref{-256,h}.
	Similarly,  \eqref{128,h} and \eqref{128,3k+2,h} are equivalent,
	and \eqref{-72,h} and \eqref{-72,3k+2,h}
	are equivalent. Also, \eqref{-25,h} and \eqref{-25,3k+2,h} are equivalent,
	and \eqref{24,h} and \eqref{24,3k+2,h} are equivalent.
	
	Now it remains to prove \eqref{-256,h}, \eqref{128,h}, \eqref{-72,h}, \eqref{-25,h} and \eqref{24,h}.
	As this can be done in the way we prove Theorem \ref{1}, we omit the details.
	\qed
	
	\medskip
	\noindent{\bf Declaration of Interests}. There are no competing interests to declare.

\end{document}